\documentstyle[12pt,cd]{article}

\newtheorem{thm}{Theorem}
\newtheorem{prop}{Proposition}

\newtheorem{defn}{Definition}
\newtheorem{lemma}{Lemma}

\newcommand{\T}{Teichm\"{u}ller }

\newcommand{\C}{$\hat {\Bbb C}$}

\newcommand{\Bbb}{\bf} 
\newcommand{\sn}{\smallskip\noindent}

\newcommand{\cdr}{\cd\rightarrow}
\newcommand{\cdl}{\cd\leftarrow}

\newcommand{\cdd}{\cd\downarrow}

\newenvironment{pf}{\noindent {\it Proof.~}}{\hfill$\Box$\newline}
\begin{document}

\title{On \T spaces of Koebe groups}
\author{Pablo Ar\'{e}s Gastesi}
\date{January 27, 1995}

\maketitle

\begin{abstract}
In this paper we parametrize the \T spaces of constructible
Koebe groups, that is
Kleinian group that arise as covering of $2-$orbifolds determined by
certain normal subgroups of their fundamental groups. We also study
the covering spaces of the \T spaces of those Koebe groups.
Finally we
prove an isomorphisms theorem similar to the Bers-Greenberg theorem
for Fuchsian groups. Our method yields a technique to compute
explictly generators of Koebe groups, possibly by programming a
computer.
\end{abstract}

\section{Introduction}
In this paper we will extend the resuls of I. Kra \cite{kra:horoc} and the
author \cite{ares:horoc1} to the class of constructible Koebe
groups. Our main goal is to produce a set of coordinates
for the \T spaces of Kleinian groups that allow explicit
computations. We also prove some results concerning the deformation
space of Koebe groups.

One of the most interesting object associated to a Riemann surface
$S$ is its moduli or Riemann space, $R(S)$, the space that
parametrizes the different complex structures on $S$ modulo
biholomorphisms.
One possible way of studying $R(S)$ is by passing to its
universal covering space, known as \T space, $T(S)$, which is
the set of complex structures on $S$ modulo isometries (we will assume that
$S$ has a metric of constant curvature $-1$)
isotopic to the identity. In order to get
explicit coordinates on $T(S)$, we can study the Riemann surface $S$ via its
universal
covering space, the upper half plane plane, obtaining in this way that the
group of deck transformations becomes a group of M\"{o}bius
transformations. The points in $T(S)$ are then given by different
subgroups of $PSL(2,\Bbb R)$ (with a preferred set of generators).
But groups acting on the upper half plane
are quite difficult to handle for computations,
so we need to represent $S$ in a slightly
different way. B. Maskit \cite{mas:class}
proved that it is possible to find a finitely
generated Kleinian group $\Gamma$, with a simply connected
invariant component $\Delta$
such that $\Delta/\Gamma\cong S$. The set of groups quasiconformally
conjugated to
$\Gamma$, modulo conjugation by M\"{o}bius transfomations, is the \T space
of $\Gamma$, $T(\Gamma)$. This space is the cartesian product of the \T
spaces of all surfaces represented by $\Gamma$. But by the same result of
B. Maskit, we have that the surfaces uniformized by $\Gamma$ other than
$S$, are rigid, i.e. their \T spaces are points. Therefore $T(\Gamma)$
becomes a model for $T(S)$ with the advantage that these groups are
good for explicit computations.
In this line lies the work of I. Kra
(\cite{kra:horoc}), where he studied the case of compact surfaces with
finitely many punctures. Later, the author extended the results to the
case of $2$-orbifolds (topological surfaces, with a complex structure,
where each point has a
neighborhood modeled over a Euclidean disc quotiented a finite group of
rotations) in \cite{ares:horoc1}. Observe that orbifolds in particular
include the Riemann surfaces, if we take the rotation group to be
trivial. In this paper we will study another
type of groups related to planar covering of $2$-orbifolds.
We start by choosing a maximal
partition $\cal P$ on the orbifold $S$; that is, a set of curves that
splits $S$ into spheres with three marked points or holes. To each curve
of the partition $a_j$, we assing an integer number (bigger than $2$) or
$\infty$, $\mu_j$, called a weight. We then consider the normal subgroup
$H$ of $\pi_1(S)$
generated by the curves $a_j^{\mu_j}$ with finite $\mu_j$. By the same
result of B. Maskit we have that the
covering determined by $H$ produces a Kleinian group $\Gamma$,
known as Koebe group,
with an invariant component $\Delta$, which is simply connected if and
only if all $\mu_j=\infty$.
The group $\Gamma$ uniformizes $S$ and
rigid orbifolds, so the study of $T(\Gamma)$ is somehow equivalent to the
study of $T(S)$. The main result in this line is the following:
\begin{thm} Let $S$ be an orbifold with signature
$(p,n;\nu_1,\ldots,\nu_n)$ satisfying $2p-2+n>0$ and $3p-3+n-\sum
1/\nu_j>0$. Let ${\cal P}=\{a_1,\ldots,a_{3p-3+n}\}$ be a maximal partition on
$S$, and let
${\cal N}=\{\mu_1,\ldots,\mu_{3p-3+n}\}$ be a set of weights.
Assume that $\Gamma$ is a Koebe group uniformizing
$(S,\cal P,\cal N)$. Then there exists a set of
coordinates,$(\alpha_1,\ldots,\alpha_{3p-3+n})$,
on the \T space $T(\Gamma)$, and a set of positive numbers,
$(r_1,\ldots,r_{3p-3+n})$, such that
$$\prod_{j=1}^{3p-3+n}U_j\subset T(\Gamma),$$
where $U_j=\{\alpha_j;~Im(\alpha_j)>r_j\}$ if $\mu_j=\infty$ or
$U_j=\{\alpha_j;~0<|\alpha_j|<r_j\}$ if $\mu_j<\infty$.
Moreover, we also have the inclusions
$$T(\Gamma)\subset\prod_{j=1}^{3p-3+n}V_j,$$
where $V_j$ is the upper half plane if $\mu_j=\infty$, or the unit
disc if $\mu_j<\infty$.
The entries of a set of generators correspoding to a point in $T(\Gamma)$
can be computed {\bf explicitly} in terms of the coordinates, and vice
versa. The numbers $r_j$ depend only on the topology of
$S,\cal P,\cal N$.
\end{thm}
A Riemann surface $S$, or more generally an orbifold, can be constructed
from more basic orbifolds, $S_1$ and $S_2$, by removing discs and
identifying their
boundaries. In the case that such identification is given by a
formula of the type $zw=t$, where $z$ and $w$ are local coordinates
on $S_1$ and $S_2$ respectively, we say that $S$ has been constructed
from $S_1$ and $S_2$ by plumbing techniques, with parameter $t$. The
above study of Koebe groups allows us to understand our computations
in terms of plumbing constructions.
\begin{thm}The orbifold corresponding to the point
$(\alpha_1,\ldots,\alpha_{3p-3+n})\in T(\Gamma)$ can be constructed
by plumbing techniques with parameters $(\tau_{1},\ldots,\tau_{3p-3+n})$,
where $\tau_j=exp(\pi ic_j\alpha_j)$ if $\mu_j=\infty$, or
$\tau_j=c_j\alpha_j$, if $\mu_j<\infty$. The constants $c_j$ depend
only on the topology of the orbifold and the partition.\end{thm}

If one of the $\mu_j's$ is finite, then $T(\Gamma)$
will not be $T(S)$, although it can be proven that the latter space is
the universal covering of the former. Our result identifies the
covering group of $T(S)\rightarrow T(\Gamma)$.

\begin{thm} The covering group of the mapping $T(S)\rightarrow T(\Gamma)$
is the normal subgroup of the mapping class group of $S$ generated by
the Dehn twists $\tau_j^{\mu_j}$ around the curves $a_j$ of the partition
with finite $\mu_j$.
\end{thm}

The Bers-Greenberg theorem tells us that the complex structure of
$T(\Gamma)$ depends only on the pair $(p,n)$ of $S$. This result was
first proven for the case of Fuchsian groups in \cite{bg:isom}.
The author in \cite{ares:bers} extended to the case of Koebe groups where all
the $\mu_j=\infty$, following \cite{ek:hol}. In this paper we generalized
it to the case of constructible Koebe groups with finite weights.
\begin{thm}Let $S$, $\cal P$ and $\cal N$ be as in theorem $1$.
Let $S_0$ be the
surface obtained by removing from $S$ all the points with finite
ramification number. Assume that $\Gamma$ and $\Gamma_0$ are two
Koebe groups uniformizing $(S,\cal P,\cal N)$ and $(S_0,\cal P,\cal N)$
respectively. Then the deformation spaces $T(\Gamma)$ and
$T(\Gamma_0)$ are conformally equivalent.\end{thm}

This paper is organized as follows: in \S $2$ we give the necessary background
on \T spaces and Kleinian groups. The Koebe groups we study are constructed
from more basic groups, known as triangle groups, which are studied in
detail in \S $3$. In \S $4$ we compute the coordinates of theorem $1$; by a
theorem of B. Maskit \cite{mas:moduli}, it suffices to consider the
cases of dimension $1$, which we do in that section, and we also indicate
how to proceed in the general situation. We will also give a relation
between the coordinates on deformation spaces and plumbing constructions
on $2$-orbifolds,
proving theorem $2$. In \S $5$ we prove theorems $3$ and $4$.

{\bf Acknowledgements.} I would like to thank Irwin Kra for introducing me in
the study of Kleinian groups and their deformation spaces. This paper is
an extension of the author's Ph. D. thesis. The main work was done while I
was visiting the University of Joensuu and the Tata Institute. I would like
to thank both places for financial support. The
results in \S\S $3$ and $4$ were obtained simultaneously and
independently by Jouni Parkkonen in \cite{jouni:gluing}.

\section{Background on Kleinian groups and \T spaces}
{\bf 2.1.} A {\bf Kleinian group} $\Gamma$ is a discrete subgroup of
$PSL(2,\Bbb C)$ such that the set of points of $\hat{\Bbb C}={\Bbb
C}\cup\{\infty\}$ where $\Gamma$ acts discontinuously, the {\bf
regular set} $\Omega=\Omega(\Gamma)$, is not empty. If \C$-\Omega$
consists of at most $2$ points, we say that the group is {\bf elementary}.

Let $\Delta$ be a invariant component of $\Gamma$; that is, a
component of $\Omega$ such that $\gamma(\Delta)=\Delta$ for all
$\gamma\in\Gamma$. If $\Gamma$ is finitely generated, the natural
mapping $\pi:\Delta\rightarrow\Delta/\Gamma$ is a covering
of a compact surface of
genus $p$ with finitely many punctures, ramified over finitely
many points. We say that $S=\Delta/\Gamma$ is an {\bf orbifold} of
signature $(p,n;\nu_1,\ldots,\nu_n)$, where the $\nu_j\in{\Bbb
Z}^+\cup\{\infty\}$, $\nu_j\geq 2$. The $\nu_j'$s are called the
ramification values. We have that $\pi$ is $\nu_j-$to$-1$ in a
neighborhood of a point $x_j$ whose ramification value is $\nu_j$. We
will assign ramification value equal to $\infty$ to the punctures of
$S$. The pair $(p,n)$ is called the {\bf type} of the orbifold. All
orbifolds in this paper, except for those of type $(0,3)$,
will satisfy the following two conditions:
\begin{equation}\left\{\begin{array}{l}
2p-2+n>0\\
3p-3+n-\sum 1/\nu_j>0\end{array}\right .\end{equation}

If $\Delta/\Gamma$ is an orbifold of signature $(0,3;\nu_1,\nu_2,\nu_3)$,
then $\Gamma$ is
called a {\bf triangle group}. Triangle groups are divided into
{\bf hyperbolic, parabolic} or {\bf elliptic}, depending on whether
$1-(1/\nu_1+1/\nu_2+1/\nu_3)$ is positive, zero or negative, respectively.
If the group is elliptic or parabolic, then it is elementary, and
$\Delta=\Omega$. If $\Gamma$ is hyperbolic, $\Omega$ consists of two
discs or half planes, and $\Delta$ is any of them.

A {\bf constructible Koebe group} is a Kleinian group $\Gamma$, with
an invariant component $\Delta$ such that $\Gamma$ can be built up
from elementary and hyperbolic triangle groups by finitely many
applications of the Klein-Maskit Combination Theorems (see
\cite{mas:klein4} for the latest version of these theorems). In
particular we get that $\Gamma$ is finitely generated. by the
signature of $\Gamma$ we will understand the signature of $\Delta/\Gamma$.
For the rest of this paper, $\Gamma$ will be a constructible Koebe group,
unless otherwise stated.

\sn {\bf 2.2.} A {\bf maximal partition with weights} on an orbifold
$S$, is a pair $(\cal P,\cal N)$ where:

(i) ${\cal P}=\{a_1,\ldots,a_{3p-3+n}\}$ is
a set of $3p-3+n$ simple closed disjoint unoriented curves on
$S_0=S-\{x_j;~\nu_j<\infty\}$,
such that no curve
of $\cal P$ bounds a disc or a punctured disc on $S_0$,
and no two curves of $\cal P$ bound a cilinder on $S_0$;

(ii) ${\cal N}=(\mu_1,\ldots,\mu_n)\in ({\Bbb Z}\cup\{\infty\})^{3p-3+n}$,
with $\mu_j\geq 3$;

(iii) the weight $\mu_j$ is assigned to the curve $a_j$.
\begin{thm}[Maskit Existence Theorem] Given an orbifold $S$ with
signature satisfying ($1$), and maximal partition $(\cal
P,\cal N)$, there exists a unique (up to conjugation in $PSL(2,\Bbb
C)$) constructible Koebe group $\Gamma$, with invariant component
$\Delta$, such that:

(i) $S\cong\Delta/\Gamma$;

(ii) to each curve $a_j$ of $\cal P$ corresponds a unique conjugacy class
of elements of $\Gamma$ generated by a transformation of order $\mu_j$;

(iii) $(\Omega-\Delta)/\Gamma$ is the union of orbifolds of type
$(0,3)$ obtained by squeezing each curve $a_j$ of $\cal P$ to a point
of ramification $\nu_j$,
and discarding all orbifolds of parabolic or elliptic signature that appear.
\end{thm}
\begin{pf}It is easy to see that the condition of theorem $X.D.15$ in
\cite[pg. 281]{mas:kg} are satisfied. Existence is then given in
theorem $X.F.1$ of the same reference, while uniqueness is theorem
$1$ of \cite{mas:class}.\end{pf}
We will say that $\Gamma$ uniformizes the triple $(S,\cal P,\cal N)$.

\sn {\bf 2.3.} Let $G$ be a finitely generated non-elementary
Kleinian group. The {\bf \T or deformation space} of $G$ is the set
$$T(G)=\{w:\hat{\Bbb C}\rightarrow\hat{\Bbb C};~w~ \mbox{is
quasiconformal},~wG w^{-1}\leq PSL(2,{\Bbb C})\}/\sim,$$
where $w_1\sim w_2$ if there is a M\"{o}bius transformation $A$ such that
$w_1 g w_1^{-1}=Aw_2 g w_2^{-1}A$, for all $g\in G$.

\noindent If $\Gamma$ is a constructible Koebe group, then
$T(\Gamma)$ is a complex manifold of dimension $3p-3+n$. In the case
that all the weights of $\Gamma$ are equal to $\infty$, we have that
$T(\Gamma)$ is equivalent to $T(S)$, the deformation space of
$S=\Delta/\Gamma$,
which is the set of quasiconformal homeomorphisms of $S$ modulo those
isotopic to conformal mappings (see \cite{kra:spaces} and
\cite{nag:teic}).

\sn {\bf 2.4.} There is a way of decomposing $\Gamma$ into simpler
groups as follows. Let $T_j$ be the connected component of
$S-\{a_k;~a_k\in {\cal P},~k\neq j\}$ containing $a_j$. Let $D_j$ be
a connected component of $\pi^{-1}(T_j)$, where
$\pi:\Delta\rightarrow S$ is the natural projection. Denote by
$\Gamma_j$ the stabilizer of $D_j$ in $\Gamma$; that is,
$\Gamma_j=\{\gamma\in\Gamma; \gamma(D_j)=D_j\}$. These groups
are Koebe group of type $(0,4)$ or $(1,1)$. Therefore
$dimT(\Gamma_j)=1$, and it is clear that there is a mapping from
$T(\Gamma)$ into $T(\Gamma_j)$ given by restriction. We can choose
these subgroups so that $\Gamma_j\cap\Gamma_{j+1}=F_j$ is a triangle
group for $1\leq j\leq 3p-4+n$.
\begin{thm}\label{embed}[Maskit Embedding,\cite{mas:moduli},
\cite{kra:variational}]
The mapping given by restriction
$T(\Gamma)\longrightarrow$\newline
\noindent$\displaystyle{\prod_{j=1}^{3p-3+n}T(\Gamma_j)}$
is holomorphic, injective and with open image as long as the groups
$F_j$ do not have signature $(0,3;2,2,\nu)$ for finite $\nu$.
\end{thm}
\begin{pf}
See \cite{kra:variational} and observe that the only triangle groups with
non trivial centralizer in $PSL(2,\Bbb C)$ have the above metioned
signatures.\end{pf}
In our case, since we are assuming that the weights are always
strictly biggerr than $2$, we always have such an embedding.

\sn {\bf 2.5.} An element $A$ of a Kleinian group $G$ is said to be
{\bf primitive} if it has no roots in $G$; that is, if $B\in G$ and
$B^n=A$, then $n=\pm 1$.

An elliptic element $C$ of finite order $n$ is conjugated in $PSL(2,\Bbb C)$
to a rotation of the form $z\mapsto e^{2k\pi i/n}z$, with $k$ and $n$
relatively primes. If $k=\pm 1$ we say that $C$ is {\bf geometric}
(\cite[pg. 96]{mas:kg}).

\section{Parametrization of triangle groups}
{\bf 3.1.} It is a classical fact that two triangle groups with the
same signature are conjugate in PSL($2,\Bbb C$) \cite[pg. 217]{mas:kg}.
Therefore,
to determine a triangle group alll we need is its signature and three
distinct points of \C, which we will call {\bf parameters}. We are
interested in a having a technique to compute
{\it explictly} generators for these groups. This requires
some canonical choices. The main goal of this section is to develop
such techniques. Later, in \S 4, we will use these
generators to compute coordinates on \T spaces of Koebe groups.

$\Gamma(\nu_1,\nu_2,\nu_3;a,b,c)$ will denote a triangle group of
signature $(0,3;\nu_1,\nu_2,\nu_3)$ with a pair of canonical
generators, $A$ and $B$, for the parameters $a,b,c$. We will always
assume that $A$ and $B$ are primitive and geometric (if elliptic), and
$|A|=\nu_1$, $|B|=\nu_2$ and $|AB|=\nu_3$. Here $|T|$
denotes the order of a M\"{o}bius transformation, with parabolic
elements considered as elements of order equal to $\infty$.
For technical
reasons (see \S 3.2.) we will always assume that $\nu_1>2$. It will be
clear from our definitions that if $T$ is a M\"{o}bius
transformation, then $TAT^{-1}$ and $TBT^{-1}$ are canonical
generators for $\Gamma(\nu_1,\nu_2,\nu_3;T(a),T(b),T(c))$. This will
simplify the proofs of this section by taking $a=\infty$, $b=0$ and $c=1$.

\sn {\bf 3.2.} If $A$ is a generator of a triangle group, then it is
either parabolic or elliptic. In both cases, there are circles (or
more generally, closed Jordan curves), invariant under $A$. Let
$\tilde a$ be one such curve. We orient it by picking a point
$z\in\tilde a$, not fixed by $A$, and then requiring that $z$,
$A(z)$ and $A^2(z)$ follow each other in the positive orientation.
This is possible since we are assuming that the cuves of
the partition are always uniformized by elements of order strictly
bigger than $2$. The cases of groups with partitions curves
uniformized by involutions will be treated in a
forthcoming publication.\vspace{2mm}

\noindent{\large\bf Hyperbolic Groups.}\vspace{2mm}

\sn {\bf 3.3.} The results of this section are taken from
\cite{ares:horoc1}, which are generalization of those in \cite{kra:horoc}.
Given three distinct points $a,b,c$ in $\hat{\Bbb C}$, let $\Lambda$ be the
circle determined by them and oriented so that $a,b,c$, follow each other
in the positive orientation. Let $\Delta$ be the disc to the
left of $\Lambda$;
that is, the set of points $z$ with $cr(z,a,b,c)>0$, where $cr$
denotes the cross ratio of four distinct point on the Riemann sphere,
chosen (there are $6$ possible definitions of cross ratio) such that
$cr(\infty,0,1,z)=z$. Let $L$ and $L'$ be the circles passing through
$\{a,b\}$ and $\{a,c\}$ respectively, and orthogonal to $\Lambda$.
\begin{defn} Given $z_1$ an $z_2$ in $\overline{\Delta}\cap L$, we
will say that they are {\bf well ordered} with respect to $(a,b,c)$
if one of the following conditions is satisfied:\\
(i) $z_1=a$,\\
(ii) $z_2=b$,\\
(iii) $z_1\neq a$, $z_2\neq b$ and $cr(a,z_1,z_2,b)>1$.
\end{defn}
For example, if $a=\infty$, $b=0$ and $c=1$, then $\Delta$ is
the upper half plane. Two points $z_1=\lambda i$ and $z_2=\mu i$
are well ordered with respect to these parameters if $\lambda>\mu$.
Na\"{\i}vely speaking, this means that $z_1$ is closer to $\infty$ than
$z_2$.

Let $\Gamma$ be a triangle group with signature $(0,3;\nu_1,\nu_2,\nu_3)$
whose limit set ($\hat{\Bbb C}-\Omega$)
is $\Lambda$, and let $A$ and $B$ be two generators of
$\Gamma$.
\begin{defn}$A$ and $B$ are the {\bf canonical generators} of
$\Gamma(\nu_1,\nu_2,\nu_3;a,b,c)$ if the following conditions are satisfied:\\
(i) $A$ and $B$ have their fixed points in $L$ and $AB$ in $L'$;\\
(ii) if $z_1$ and $z_2$ are the fixed points of $A$ and $B$ in
$\overline{\Delta}\cap L$, then they are well ordered with respect to
$(a,b,c)$;
\end{defn}

\begin{prop}There exists a unique
$\Gamma(\nu_1,\nu_2,\nu_3;a,b,c)$.
In the case of $a=\infty$, $b=0$ and $c=1$, the canonical generators
$A$ and $B$ Are given by:

(1) $\nu_1=\infty$,
$$A=\left[ \begin{array}{cc}-1 &-2\\0&-1\end{array}\right],~
B=\left[ \begin{array}{cc} -q_2&b\\q_2+q_3&-q_2\end{array}\right],$$
where $q_i=\cos(\pi/\nu_i)$, $b=\frac{q_2^2-1}{q_2+q_3}$;

(2) $\nu_1\neq\infty$,
$$A=\left[ \begin{array}{cc}-q_1 &-kp_1\\k^{-1}p_1&-q_1\end{array}\right],~
B=\left[ \begin{array}{cc} -q_2&-hp_2\\h^{-1}p_2&-q_2\end{array}\right],$$
$q_i=\cos(\pi/\nu_i)$, $p_i=\sin(\pi/\nu_i)$,
$$k=\frac{q_2+q_1q_3+q_1l}{p_1l},~
h=\frac{kp_1p_2}{q_1q_2+q_3+l}, l=\sqrt{q_1^2+q_2^2+q_3^2+2q_1q_2q_3-1}>0.$$
\end{prop}
{\bf Remarks.} 1. Observe that the generators
$A$ and $B$ in the second case converge to
those of the first case as $\nu_1\rightarrow\infty$.

2.- The proof of this result (see \cite{ares:thesis}) is based on
a detalied study of hyperbolic triangles in the upper half
plane.

3.- Wrom the above expression we can see that $A$ ($B$) fixes the points
$\pm ki$ ($\pm hi$). There is a
relation between these formul\ae~ and the geometry of the orbifold
$S={\Bbb H}/\Gamma$ as follows. Put on $S$ the natural metric of
constant curvature $-1$ coming the Poincar\'{e} metric on the upper
half place. Then the distance from $P_1$ to $P_2$ is
$\log(\frac{k}{h})$, where $P_j$ is the point with ramification
$\nu_j$, $j=1,2$, and we are assuming that $\nu_j<\infty$.\vspace{2mm}

\noindent{\large\bf Parabolic Groups.}\vspace{2mm}

\sn {\bf 3.4.} Let us consider first the case of signature
$(0,3;\infty,2,2)$.
\begin{defn}$A$ and $B$ are {\bf canonical generators}
of $\Gamma(\infty,2,2;a,b,c)$,
if the folowing conditions are satisfied:\\
$A(z)=a$, $B(b)=b$ and $AB(c)=c$;
\end{defn}
\begin{prop}
There exists a unique $\Gamma(\infty,2,2;a,b,c)$
\end{prop}
\begin{pf} Conjugate to get $a=\infty$, $b=0$ and $c=1$.
Then $A(z)=z+\alpha$,
for some complex number $\alpha$. Since all elements of the group must
fixed $\infty$ (\cite[pg. 91]{mas:kg}), we get that $B(z)=-z$, forcing
$\alpha$ to be equal to $2$.\end{pf}

\sn {\bf 3.4.} For the other parabolic signatures, as well as for
the ellliptic cases, we need a new concept due to J. Parkkonen
\cite{jouni:gluing}.
\begin{defn}Let $M$ be a M\"{o}bius transformation of finite order strictly
bigger than $2$, and let $x$ be a fixed point of $M$. Let $z$ be any point
not fixed by $M$. We will say that $x$ is the {\bf right} ({\bf left})
fixed point of $M$ if the cross ratio $cr(x,z,M(z),M^2(z))$ has positive
(negative) imaginary part.\end{defn}
\begin{lemma}
The above definition is invariant under conjugation by elements of
$PSL(2,\Bbb C)$.
\end{lemma}
This lemma simply means that if $x$ is the right fixed point of $M$, then
$T(x)$ is the right fixed point of $TMT^{-1}$ for any M\"{o}bius
transformations $T$.
\begin{lemma}The above definition does not depend on the point $z$.
\end{lemma}
\begin{pf} By the previous result we can assume that $M(z)=kz$, with $|k|=1$.
Then an easy computation shows that
$cr(\infty,z,kz,k^2z)=k+1$ while
$cr(0,z,kz,k^2z)=1+\frac{1}{k}$, for all $z\neq \infty,0$.\end{pf}

\sn
{\bf 3.6.} We are now in a position to define the canonical generators
for the rest of the parabolic signatures.
\begin{defn}The {\bf canonical generators} $A$ and $B$ of
$\Gamma(\nu_1,\nu_2,\nu_3;a,b,c)$, where all the ramification values are
finite, satisfy the following conditions:\\
(i) $a$ is the right fixed point of $A$ and $b$ is its left fixed point;\\
(ii) $B$ fixes $a$ and $c$;
\end{defn}
\begin{prop}There exists a unique $\Gamma(\nu_1,\nu_2,\nu_3;a,b,c)$, with
the signature in the conditions of the above definition.
\end{prop}
\begin{pf} $A(z)=\lambda z$, with $\lambda=exp(2\pi i/\nu_1)$, assuming
that the parameters are $\infty,0,1$ as usual. $B$ must be of the form
$B(z)=\mu(z-1)+z$, where $\mu=exp(\pm 2\pi i/\nu_2)$. Using the
fact that $\frac{1}{\nu_1}+\frac{1}{\nu_2}+\frac{1}{\nu_3}=1$, one can
easily that $\mu=exp(2\pi i/\nu_2)$, since
otherwise $AB$ dose not have the correct order.\end{pf}\vspace{2mm}

\noindent{\large\bf Elliptic Groups}\vspace{2mm}

\sn {\bf 3.7.} By a result of I. Kra and B. Maskit \cite{km:defor},
the elliptic triangle groups of signature $(0,3;2,2,\nu)$, with finite $\nu$,
cannot be parametrized. Nevertheless, we can have a definition of canonical
generators, although not a uniqueness statement.
\begin{defn}The {\bf canonical generators} $A$ and $B$ of
$\Gamma(\nu,2,2;a,b,c)$ satisfy:\\
(i) $A$ fixes $a$ and $b$, and $a$ is the right fixed point;\\
(ii) $B(c)=c$;
\end{defn}
Since $B$ has order $2$, there is no way to differentiate between its two
fixed points. This is why we do not have a uniqueness statement. For
computational purposes, the following result is enough, although it does not
guarantee coordinates on \T space (see remark 3 in \S 4.4).
\begin{prop}Given three distinct points $a,b,c$ on the Riemann sphere, there
is a unique point $d\in\hat{\Bbb C}-\{a,b,c\}$ such that
$\Gamma(\nu,2,2;a,b,c)=\Gamma(\nu,2,2;a,b,d)$.
\end{prop}
\begin{pf} We first compute to get $A(z)=e^{2\pi i/\nu}z$ and $B(z)=1/z$, if
$a=\infty$, $b=0$, and $c=1$. Since $B$ also fixes $-1$, we get that
$d=-1$.\end{pf}

\sn {\bf 3.8.} The last cases are those elliptic groups with exactly one
point of ramification values $2$.
 Without of generality, we can ordered the
ramification values so that the signature is of the form
$(0,3;\nu_1,\nu_2,2)$. This
avoids some long computations and does not loose any mathematical insight.
\begin{defn}The {\bf canonical generators} of
$\Gamma(\nu_1,\nu_2,2;a,b,c)$ satisfy:\\
(i) $A$ fixes $a$ and $b$, and $a$ is the right fixed point;\\
(ii) the right fixed point of $B$ is $c$;
\end{defn}
\begin{prop}
There is a unique $\Gamma(\nu_1,\nu_2,2;a,b,c)$ in the above conditions.
\end{prop}
\begin{pf} As usual, we make $a=\infty$, $b=0$ and $c=1$ by conjugation.
Let $q_j$ and $p_j$ denotes $\cos(\pi/\nu_j)$ and $\sin(\pi/\nu_j)$,
$j=1,2$, respectively. Then
$A(z)=\lambda^2z$, where $\lambda=exp(\pi i/\nu_1)$. Since $AB$ has
zero trace, we can assume that $B$ has negative trace. If $B$ has a matrix
expression given by $\left(\begin{array}{cc}a&b\\c&d\end{array}\right)$,
then we have:
$$\left\{\begin{array}{rcl}a+b&=&c+d\\ad-bc&=&1\\a+d&=&-2q_2\\
-\lambda a-\lambda^{-1}d&=&0\end{array}\right .$$
The last two equations give
$\displaystyle{a=\frac{-q_2\overline{\lambda}i}{p_1}}$.
We get
$b=-q_2-a\pm ip_2$ and $c=q_2+a\pm ip_2$.
A simple computation shows that
$$cr(1,\infty,\frac{a}{c},\frac{a^2+bc}{ac+dc})=
1+\frac{1}{-1+2q_2^2\mp 2ip_2q_2}.$$
Since $1$ has to be the right fixed point of $B$, we get that
$b=-q-a-ip$ and $c=q+a-ip$.\end{pf}\vspace{2mm}

\noindent{\large\bf Hyperbolic Groups Revisited}\vspace{2mm}

\sn {\bf 3.9.} One can re-write the results about canonical generators of
hyperbolic groups using the concept of right and left fixed points.
The existence and uniqueness proposition is then as follows:
\begin{prop}
Let $(0,3;\nu_1,\nu_2,\nu_3)$ be a hyperbolic signature where all the
ramification values are finite. Then $A$ and $B$ are the canonical generators
of $\Gamma(\nu_1,\nu_2,\nu_3;a,b,c)$ if and only if:\\
(i) $A$ and $B$ have their fixed points in $L$, and $AB$ has them in $L'$;\\
(ii) the left fixed point of $A$ lies in $\Delta$;
\end{prop}
The proof is a long but easy computation left to the reader.\vspace{2mm}

\noindent{\large\bf Local Coordinates}\vspace{2mm}

\sn {\bf 3.10.} In order to relate the parametrization of triangle
groups and the geometry of orbifolds (\S 4.9), we need to introduce
some local coordinates on the basic orbifolds of type $(0,3)$. These
coordinates were first found by I. Kra in \cite{kra:horoc}. There he
considered an orbifold $S$ of signature $(0,3;\infty,\infty,\infty)$,
with punctures $P_1$, $P_2$ and $P_3$. $S$ has a metric of constant
curvature $-1$ coming from its universal covering space, the upper
half plane $\Bbb H$. As the uniformizing group we can take
$\Gamma(\infty,\infty,\infty;\infty,0,1)$, generated by $A(z)=z+2$ and
$B(z)=-z/(2z-1)$. Then, if $\xi\in\Bbb H$, the function
$z(\xi)=exp(\pi i\xi)$ is invariant under $A$ and induces a local
coordinate on a neighborhood of $P_1$ in $S$. We can consider for
example those $\xi$ with imaginary part bigger than $1$. In $S$ there
is a unique geodesic $c$ such that, when it is parametrized by the
arc length, $c$ satisfies $lim_{s\rightarrow +\infty}c(s)=P_1$ and
$lim_{s\rightarrow -\infty}c(s)=P_2$. $z$ is charaterized by mapping
$c$ isometrically into the unit interval $(0,1)$ in the punctured
disc (with its natural Poincar\'e metric). We say that $z$ (or more
precisely, the germ of holomorphic functions determined by it) is a {\bf
preferred coordinate} at $P_1$ relative to $P_2$.

\noindent In the case of
hyperbolic groups with torsion, we still have uniqueness of geodesic
and coordinates. At a point $P_1$ of finite ramification value
$\nu_1$, the preferred coordinate looks like $z(\xi)=\xi^{\nu_1}$,
where $\xi$ lies in a neighborhood of $0$ in $\Bbb C$, and the group
we are considering is $\Gamma(\nu_1,\nu_2,\nu_3;-1,1,\frac{-1+ki}{1+ki})$.

\noindent In the parabolic and elliptic cases, we do not have uniqueness
statements for preferred coordinates, but nevertheless, we can define
coordinates by giving a function that generates the correspoding germ
of holomorphic mappings. For the parabolic group
$\Gamma(\infty,2,2;\infty,0,1)$, we define the coordinate at the
puncture $P_1$ relative to $P_2$ (one of whose lifts to $\Bbb C$ is
$0$) by $z(\xi)=exp(\pi i\xi)$. For parabolic and elliptic
groups $\Gamma(\nu_1,\nu_2,\nu_3)$ with $\nu_1<\infty$, we have to
take a more general concept of coordinates by allowing our local
patch to be modelled on the Riemann sphere. By this we mean that a
preferred local coordinate on an orbifold will mapped a neighborhood
of the special point into a disc centered at $0$,
with the point being sent to the
origin, or into the exterior of a disc (also centered at $0$), with the
image of the point under consideration being the point $\infty$.
We then have that the preferred coordinates at $P_1$ and relative to
$P_2$ are given by either
$z(\xi)=\xi^\mu$, or $z(\xi)=(1/\xi)^\mu$.

As the name suggests, for each of the groups of this
section we can construct a fundamental domain for their action on the
regular set (if the group is elementary) or in one of the two
components of $\Omega$ (in the case of hyperbolic groups) by taking a
triangle with angles $\pi/\nu_j$, and then reflecting it in one of
its sides. For example, in the case of the group
$\Gamma(4,4,2;\infty,0,1)$ we get a rectangle with vertices at the
points $0$, $(1+i)/2$, $1$ and $(1-i)/2$. Consider now the point $0$,
which project to a point with ramification number $4$, say $P_1$.
Then we can see that the distance from $0$ to the line joining
$1$ and $(1+i)/2$ is $1$. This means that the disc $\{z:~|z|<1\}$
maps into a neighborhood of $P_1$ on the quotient orbifold. Similarly, from
we have discs of radius $1$ around the other two ramification points.
In general, we have the following result:
\begin{lemma}\label{radius}
Let $\Gamma(\nu_1,\nu_2,\nu_3;a,b,c)$ be a triangle group. Let $P_j$
be a ramification point in the in the quotient orbifold,
$j=1,2,3$. Let $p_j$ be a lifting of $P_j$. Then we can find a
positive number $r$ such that the disc of radius $r$ around $p_j$
projects onto a neighborhood of $P_j$. If $P_j$ is a puncture, then
we can find a horodisc around $p_j$ (i.e., a disc such that $p_j$
lies in its boundary) that projects onto a punctured disc on $S$
containing $P_j$.
\end{lemma}

\section{Coordinates on \T spaces of Koebe groups}
{\bf 4.1} In this section we will explain how to construct the Koebe
groups given by the Maskit Existence Theorem. We will also give global
coordinates for the deformation spaces of these groups, and explain the
relation between our coordiantes and
plumbing parameters. The computations of \S $3$ allows us to
construct an algorithm from which one can get {\bf explicitly}
generators for Koebe groups. This technique was used for I. Kra
\cite{kra:horoc} and the author \cite{ares:horoc1} to compute
formul\ae~ for isomorphisms between \T space. See also \cite{jouni:gluing}.

\noindent By the Maskit Embedding Theorem, the \T space of a Koebe
group can be embedded into the product of one dimensional space.
These latter sets correspond to the groups of type $(04,)$ and
$(1,1)$, which we will work out in detail, and then indicate how to
treat the general case.\vspace{2mm}

\noindent{\large\bf The $(0,4)$ case}\vspace{2mm}

\sn {\bf 4.2} Let $S$ be an orbifold of signature
$(0,4;\nu_1,\ldots,\nu_4)$. A maximal partition on $S$ consists on a
curve $a_1$, with weight $\mu$. Let $S_1$ and $S_2$ be the two parts
of $S-\{a_1\}$. Orient $a_1$ such that $S_1$ lies to its right. If we
cut $S$ along the partition curve, and glue to each resulting
boundary a disc whose center is a point with ramification value
$\mu$, then we obtain that $S_1$ and $S_2$ have been completed to
orbifolds of signatures $(0,3;\mu,\nu_1,\nu_2)$ and
$(0,3;\mu,\nu_3,\nu_4)$ respectively. Let $\Gamma_i$ be
triangle groups uniformizing $S_i$, $i=1,2$. To recover $S$ be have to do the
opposite construction: first we must remove discs from $S_1$ and
$S_2$ and then glue along the boundaries. This implies that the
elements uniformizing $a_1$ in $S_1$ and $S_2$ must be the same
M\"{o}bius transformation. In other words,
$\Gamma_1\cap\Gamma_2=<A>$, where $A$ is primitive in both groups.
The First Combination Theorem \cite[VII.C.2, pg 149]{mas:kg} tells us
that if we choose the triangle groups properly, then the group
$\Gamma=\Gamma_1*_{<A>}\Gamma_2:=<\Gamma_1,\Gamma_2>$ is a Koebe
group of the desired signature. By the classical theory of
quasoncformal mapping we have that any orbifold
of type $(0,4)$ can be uniformized by this method.

\sn {\bf 4.3.} We will work out two examples of the above
construction, one with weight equal to $\infty$ and the other with finite
weight. Let us start by the former case. Assume that $\Gamma_1$ has
hyperbolic signature $(0,3;\infty,\nu_1,\nu_2)$ and $\Gamma$ is a
infinite dihedral group, with signature $(0,3;\infty,2,2)$. We can
start with $\Gamma_1=\Gamma(\infty,\nu_1,\nu_2;\infty,0,1)$. Since
$S_2$ lies to the right of the partition curve we have that
$\Gamma_2$ must have $A^{-1}$ as one of its generators. This implies
that $\Gamma_2$ has to be conjugate to
$\Gamma(\infty,2,2;0,\infty,1)$ by a transformation $T$ such that
$TAT^{-1}=A^{-1}$. Therefore we get $T(z)=z+\alpha$. The fact that
$S_1$ lies to the right of $S_1$ implies that ${\rm Im}(\alpha)>0$. We get
that the Koebe groups uniformizing orbifolds of the above signature
with $\infty$ weight are given by the AFP construction
$$\Gamma=\Gamma(\infty,\nu_1,\nu_2;\infty,0,1)*_{<A>}
\Gamma(\infty,2,2;0,\infty,\alpha),$$  for a proper
choice of $\alpha$. Actually, if ${\rm Im}(\alpha)>1$, we can take the
line $\{z;~{\rm Im}(z)=\frac{1}{2}{\rm Im}(\alpha)\}$
as the invariant curve needed to apply
the First Combination Theorem, and the resulting group is a Koebe
group of the desired type.

Choose an $\alpha_0$ such that $\Gamma_0=\Gamma_{\alpha_0}$ is a
Koebe group, e.g. $\alpha_0=i$. Then it is clear that $\alpha$ is a
global coordinate on the deformation space $T(\Gamma_0)$. This result
is also given in \cite{km:defor} and \cite{kra:variational}, altough there the
computations are not explicit. $\alpha$ has a $PSL(2,\Bbb C)$
invariant expression given by $\alpha=cr(\infty,0,1,\alpha)$. This
means that if $\Gamma$ is a Koebe group of the right signature and
infinity weight, given by the AFP $\Gamma=\Gamma(\infty,\nu_1,\nu_2;a,b,c)
*_{<A>}\Gamma(\infty,2,2;d,e,f)$, then $\Gamma$ is \T equivalent to
$\Gamma_\alpha$, where $\alpha=cr(a,b,c,e)$.

\sn {\bf 4.4.} As an example of the AFP construction with finite
weight, we take the case of a hyperbolic group with
signature $(0,3;4,\nu_1,\nu_2)$ and a parabolic group of signature
$(0,3;4,4,2)$. We have that $\Gamma_1$ is conjugate to
$\Gamma(4,\nu_1,\nu_2;\infty,0,1)$. For practical purposes, we choose
the transformation $M(z)=\frac{-z+ki}{z+ki}$, where $k$ is given in
proposition 1, to do such conjugation. In this way we get that $A(z)=iz$ is a
canonical generator of
$\Gamma_1=\Gamma(4,\nu_1,\nu_2;-1,1,\frac{-1+ki}{1+ki})$. $\Gamma_2$
will be conjugate to $\Gamma(4,4,2;0,\infty,1)$ by a mapping of the
form $T(z)=\alpha z$. The orientation requirements imply that
$|\alpha|<1$. We then have that $\alpha$ is a coordinate on the
corresponding deformation space, and its invariant expression is
given by $\alpha=\frac{\beta-ki}{-\beta-ki}$, where
$\beta=cr(-1,1,\frac{-1+ki}{1+ki},\alpha)$.

We can see in this example that $T(\Gamma_0)$ is not $T(S_0)$
($S_0=\Delta_0/\Gamma_0$). While $T(S_0)$ is simply connected, it is
not hard to see that $0<|\alpha|<1$ and the circle $|\alpha|=r$ is
contained in $T(\Gamma_0)$ for small values of $r$. See below for
an explicit estimate of these values.

\sn {\bf 4.5.} The other cases of type $(0,4)$ are handled in a
similar way. Here we include a table with the results. See the remars
after it for more information.

\begin{center}
\begin{tabular}{||c|c|c|c|c|c||}	\hline\hline
$\Gamma_1$& $\Gamma_2$& weight&param. $\Gamma_1$&param. $\Gamma_2$&inv.
expression\\ \hline
hyp&hyp&$\infty$&$(\infty,0,1)$&$(\infty,\alpha,\alpha-1)$
&$\beta+1$ \\ \hline
hyp&hyp&$2<\mu<\infty$&$(-1,1,\displaystyle{\frac{-1+ki}{1+ki}})$
&$(-\alpha,\alpha,\alpha\displaystyle{\frac{1+ki}{-1+ki}})$
&$\displaystyle{\frac{(ki-1)(\beta-ki)}{-\beta(1+ki)+k^2-ki}}$\\ \hline
hyp&par&$\infty$&$(\infty,0,1)$&$(\infty,\alpha,\alpha-1)$
&$\beta+1$ \\ \hline
hyp&par&$3,4,6$&$(-1,1,\displaystyle{\frac{-1+ki}{1+ki}})$
&$(0,\infty,\alpha)$&$\displaystyle{\frac{\beta-ki}{-\beta-ki}}$ \\ \hline
hyp&ell&$3,4,5$&$(-1,1,\displaystyle{\frac{-1+ki}{1+ki}})$
&$(0,\infty,\alpha/x)$&$\displaystyle{\frac{4kxi}{\beta x+2kxi}}$ \\ \hline
par&par&$\infty$&$(\infty,0,1)$&($\infty,\alpha,\alpha-1)$&$\beta+1$\\ \hline
par&par&$3,4,6$&$(\infty,0,1)$ & $(0,\infty,\alpha)$&$\beta$ \\ \hline
par&ell&$3,4$&$(\infty,0,1)$&$(0,\infty,\alpha/x)$
&$\beta x$ \\ \hline
ell&ell&$3,4,5$&$(\infty,0,1)$&$(0,\infty,\alpha/x)$
&$\beta x$ \\ \hline\hline
\end{tabular}\end{center}\vspace{2mm}

\sn {\bf Remarks.} 1. In the above table $x=(q_1q_2-p_1p_2)/(q_1q_2+p_1p_2)$
($q_j=\cos)\pi/nu_j)$, $p_j=\sin(\pi/\nu_j)$).\\
2. If $\Gamma=\Gamma(\mu,\nu_1,\nu_2;a,b,c)*_{<D>}
\Gamma(\mu,\nu_3,\nu_4;d,e,f)$, then $\beta=cr(a,b,c,f)$
and the coordinate on the \T space is given by the above invariant
expression (last column).\\
3. If one of the triangle groups involved in the construction of $S$
has elliptic signature $(0,3;\nu,2,2)$, we can still compute the
Koebe groups by the above techniques. But since we do not have
uniqueness of parameters for these triangle groups, we do not obtain
a coordinate on the deformation spaces. We will not work any further
these cases.

If $\mu=\infty$, then the set $\{\alpha;~{\rm Im}(\alpha)>1\}$ is contained
in the deformation space of the corresponding Koebe group. For the
case of finite weights, let us consider a fundamental domain for
$\Gamma_1$ containing the origin. Let $d(\mu,\nu_1\nu_2)$ denote the
radius given in lemma 3 of \S 3.10. Similarly, let
$D(\mu,\nu_3,\nu_4)$ denote the radius of a disc centered at $0$ such
that the fundamental domain of $\Gamma_2$ is contained in that disc.
For example,
if the signature of $\Gamma_1$ is hyperbolic, one can take
$D(\mu,\nu_3,\nu_4)=1$. Then we have that the set
$\{\alpha;~|\alpha|<d(\mu,\nu_1,\nu_2)/D(\mu,\nu_3,\nu_4)\}$ is
contained in the \T space of the group under consideration.\vspace{2mm}

\noindent{\large\bf The $(1,1)$ case}

\sn {\bf 4.6.} To construct a Koebe group $\Gamma$ of signature
$(1,1;\nu)$ and weight $\mu$, we start with a triangle group
$\Gamma_1=\Gamma(\mu,\mu,\nu;a,b,c)$. We then remove two discs around
the two points of ramification value $\mu$, and glue their
boundaries. At the group level, this is reflected on finding a
transformation $C$ that conjugates $B^{\pm 1}$ to $A$.
Choose an orientation of $a_1$, the partition curve, so
that the ramification point corresponding to $A$ ($B$) lies to the
left (right) of the invariant lift $\tilde{a}_1$ of $a_1$. This
implies that the conjugation is of the form $CB^{-1}C^{-1}=A$. Maskit Second
Combination Theorem \cite[VII.E.5, pg 161]{mas:kg} tells us that a
proper choice of $C$ will produce a group
$\Gamma=\Gamma_1*_C:=<\Gamma_1,C>$ of the desired signature.

\sn {\bf 4.7.} Let us work the most complicated case of the above
construction, namely that of $\Gamma_1$ having hyperbolic signature
$(0,3;\mu,\mu,\nu)$, with finite weight $\mu$. We start with
$$\Gamma_1=\Gamma(\mu,\nu,\nu;-1,2,(-1+ki)/(1+ki))$$.
The canonical generators of this group are given by
$$A=\left[\begin{array}{cc}-q-ip&0\\0&-q+pi\end{array}\right],
B=\left[\begin{array}{cc}-q-pmi&pmi\\-pmi&-q+pmi\end{array}\right],$$
where $q=\cos(\pi/\mu)$ and $p=\sin(\pi/\mu)$, $m=k/(2h)+h/(2k)$,
$n=k/(2h)-h/(2k)$, and $k$ and $h=p/c$ are given in proposition 1 of \S
3.3. Observe that $m^2-n^2=1$. Actually, these two numbers can be
understood in terms of hyperbolic cosines and sines of some geometric
objects on the quotient orbifold, but we are not interested on this
line of thought. See \cite{jouni:gluing} for more information.It is
not hard to check that the transformation $T$ given by
$$T=i\left[\begin{array}{cc}R&\frac{(1-m)R}{n}\\
\frac{n}{2R}&\frac{-1-m}{2R}\end{array}\right],$$
where $R=\sqrt{\frac{m+1}{2}}$,
conjugates $B$ into $A$ (and vice versa). Then, the transformation
$C$ needed for the HNN extension will be of the form
$C_\tau=D_\tau T$, where $D_\tau(z)=\tau^2/z$. Or in a single
expression we have
$$C=i\left[\begin{array}{cc}\frac{\tau n}{2R}&\frac{-\tau (m+1)}{2R}\\
\frac{R}{\tau}&\frac{(1-m)R}{\tau n}\end{array}\right].$$
An application of the Second Combination Theorem shows that
$\Gamma=\Gamma_1*_C$ is a Koebe group of the desired type, for a
proper choice of $\tau$. It is clear that $\tau^2$ is a global
coordinate on the corresponding deformation space. Its invariant
expression is given by
$\displaystyle{(\frac{-n}{2})\frac{-\beta+ki}{\beta+ki}},$ with
$\beta=cr(-1,1,\frac{-1+ki}{1+ki},C(-1))$.

To give a bound on the value of $\tau^2$, re-write $C_\tau$ as
$C_\tau(z)=\tau^2D(z)$, with\\
$D=\left[\begin{array}{cc}1&(1-m)/n\\-n/1&(m-1)/2\end{array}\right].$
The circle of radius $s$ centered at the fixed poing of $B$,
$x=(m+1)/n$, is mapped onto a circle centered at the origin with
radius $|\tau|^2||D(x-s)|$. If these two circles are disjoint, the
Second Combination Theorem can be applied. Using basic calculus, we
get that the function $|x-s|/|D(x-s)|$ is decreasing as $s$ goes to
$0$, and its minimum values is $(m-1)/2$. Therefore, the set
$\{\tau;~|\tau|^2<(m-1)/2\}$ is contained in the \T space that we are
studying.

\sn {4.8.} As in the previous case, we include the results of the
case $(1,1)$ in the following table, where $\beta$, $k$, $m$ and
$n$ are as in \S 4.7, and
$q=\cos(\pi/\mu)$.

\begin{center}
\begin{tabular}{||c|c|c|c||}	\hline\hline
$\Gamma_1$&weight&$C$&inv. expression\\ \hline
$\Gamma(\infty,\infty,\nu;\infty,0,1)$&$\infty$&
$i\left[\begin{array}{cc}\tau&\sqrt{2/(1+q)}\\ \sqrt{(1+q)/2}&
0\end{array}\right]$&$\displaystyle{\sqrt{\frac{1+q}{2}}}\beta$ \\ \hline
$\Gamma(\mu,\mu,\nu;-1,1,\displaystyle{\frac{-1+ki}{1+ki}})$&$2<\mu<\infty$
&$i\left[\begin{array}{cc}\tau&-(m+1)\tau/n\\-n/(2\tau)&
(m-1)/(2\tau)\end{array}\right]$
&$(\displaystyle{\frac{-n}{2})\frac{-\beta+ki}{\beta+ki}}$\\ \hline
$\Gamma(4,4,2;\infty,0,1)$&$4$
&$i\left[\begin{array}{cc}0&\tau\\-1/\tau&1/\tau\end{array}\right]$
&$\beta$ \\ \hline
$\Gamma(3,3,3;\infty,0,1)$&$3$
&$i\left[\begin{array}{cc}0&\tau\\-1/\tau&1/\tau\end{array}\right]$
&$\beta$ \\ \hline
$\Gamma(3,3,2;\infty,0,1)$&$3$&
$i\left[ \begin{array}{cc}\tau&-\tau\\-2/(3\tau)&-1/(3\tau)\end{array}\right]$
&$-2/(3\beta)$\\ \hline\hline
\end{tabular}\end{center}\vspace{2mm}
The coordinate on the deformation space is $\tau$ in the first case
and $\tau^2$ in the other cases. It can be obtained by the above
invariant expression, where $\beta=cr(a,b,c,C(a))$, for a group given by
$\Gamma=\Gamma(\mu,\mu,\nu;a,b,c)*_C$.\vspace{2mm}

\noindent{\large\bf Plumbing constructions}\vspace{2mm}

\sn {\bf 4.9.} The above group theoretical constructions have a nice
geometrical interpretation as follows.
Let $S_1$ and $S_2$ be two orbifolds, not necessarily
distinct, of type $(0,3)$, and
let $x_1$ and $x_2$ be two points or punctures on $S_1$ and $S_2$
respectively, with equal ramification values.
Assume that $z_1$ and $z_2$ are local coordinates on $S_1$ and $S_2$
as given in \S 3.10. Choose a complex
number $t$ such that $\overline{U_i}$ is contained
in $S_i$, for $i=1,2$, where $U_i=\{p\in S_i;~|z_i(p)|<\sqrt{|t|}\}$.
If $S_1=S_2$, then we must also require that
$\overline{U_1}\cap \overline{U_2}=\emptyset$.
Let $S=S_1\sqcup S_2/\sim$, where
$p_1\in S_2$ and $p_2\in S_2$ are equivalent, $p_1\sim p_2$, if
$|z_i(p_i)|=\sqrt{|t|}$ and $z_1(p_1)z_2(p_2)=t$. Then $S$ is an orbifold
of type $(0,4)$ or $(1,1)$, depending on whether $S_1$ and $S_2$ are distinct
or not. Na\"{\i}vely speaking, we are removing discs from the orbifolds
$S_1$ and $S_2$ and gluing them with a twist. We say that $S$ has been
contructed from $S_1$ and $S_2$ with {\bf plumbing
parameter} $t$. See \cite{kra:horoc} for a more
detailed explanation of plumbing techniques in
the context of Riemann surfaces.

\noindent If the plumbing construction requires
two diferent orbifolds, then we take as $z$ the preferred coordinate
on $S_1$ centered at the first ramification point and relative to the
second ramification point; and similarly for $w$ on $S_2$. If we have
$S_1=S_2$, then we take as $z$ the preferred coordinate centered at
the first ramification point and relative to the second, while we
take as $w$ the preferred coordinate at the second ramification point
and relative to the first one. To compute $w$ all we have to do is find a
transformation $T$ that conjugates $B$ to $A$ and such that $T^2=id$.
This guarantees that $T$ conjugates the triangle group
$\Gamma_1$ to itself. This choice of coordinates is different to that of
\cite{kra:horoc} and \cite{ares:horoc1}, but it simplifies the
computations and agress with \cite{jouni:gluing}.
The construction with $w$ being the coordinate centered
at the second ramification point and relative to the third point gives also
a plumbing construction.

Let us compute the plumbing parameters of the
examples \S\S 4.3, 4.4 and 4.6. In the first case, we have that the
coordinate $z$ is given by $z_1=exp(\pi i\xi$, where $\xi$ lies in the
upper helf plane. Similarly, we have $z_2exp(\pi i (\alpha-\xi))$, and
therefore we get the plumbing parameter is given by $t=z_1z_2=exp(\pi
i\alpha)$. In the situation of \S 4.4, the
coordinates $z_1$ and $z_2$ are given by
$z_1(\xi_1)=\xi_1^4$ and
$z_2(\xi_2)=(\alpha/\xi_2)^4$, where $\xi_i$ is a
point in a fundamental domain of $\Gamma_i$,
$i=1,2$. The boundary identification gives $t=\alpha^4$.

For the case of tori, \S 4.6, we have
$z_1(\xi)=\xi^\mu$.
The transformation $T$ that conjugates $B$ to $A$ was found in \S 4.7.
The coordinate $w$ is then given by $w(\xi)=z(T(\xi))$. Since $C$
identifies a curve invariant under $B$ with a curve invariant under
$A$, we have that the plumbing parameter is computed by means of the
expression $t=z(C(\xi))\,w(\xi)$. Recalling the transformation $C$
from \S 4.7 we get the value $t=(\tau^2)^\mu$.

\noindent The values of the plumbing parameters for the $(0,4)$ cases are
given by $exp(\pi i\alpha)$ in the case of infinite weight, or $\alpha^\mu$
for finite weight $\mu$. In this latter case, if one of the triangle
groups is elementary and corresponds to $S_1$, then we take a preferred
coordinate that maps the ramification point to the origin in $\Bbb
C$, while if the elementary group corresponds to $S_2$, our local
coordinate around the ramification point will map such point to the
point $\infty$ in $\hat{\Bbb C}$.

\noindent The cases of tori are included in the
following table.

\begin{center}
\begin{tabular}{||c|c|c||}	\hline\hline
$\Gamma_1$&$T(\xi)$&plum. par. \\ \hline
$\Gamma(\infty,\infty,\nu)$&$\displaystyle{\frac{-2}{(1+q)\xi}}$
&$exp(\pi i\tau\sqrt{\frac{2}{1+q}})$\\ \hline
$\Gamma(\mu,\mu,\nu)$&$(\displaystyle{\frac{1+m}{n}})
\displaystyle{\frac{n\xi+1-m}{n\xi-1-m}}$
&$(\tau^2)^\mu$\\ \hline
$\Gamma(4,4,2)$&$1-\xi$&$\tau^8$ \\ \hline
$\Gamma(3,3,3)$&$1-\xi$&$\tau^6$ \\ \hline
$\Gamma(3,3,2)$&$\displaystyle{\frac{2\xi+1}{2\xi-2}}$
&$-\frac{27}{8}\tau^6$\\ \hline\hline
\end{tabular}\end{center}\vspace{2mm}

\noindent{\large\bf The general case}\vspace{2mm}

\sn {\bf 4.10.} Let us finally explain how to compute coordinates for
Koebe groups in general and vice versa; that is, given a point in the
deformation space, how to find the generators of the
corresponding group. Let $\Gamma$ be a Koebe group uniformizing an
orbifold $S$ with maximal partition with weights $(\cal P,\cal N)$.
In the classical work on deformation spaces of Fuchsian
groups, one starts with a group, say $F_0$, and then to compute the
point correspoding to any other group, $F$, one has to measure how
far $F$ is from $F_0$, via the maximal dilatation of quasiconformal
mapppings. This implies that the starting group
$F_0$ plays a special point with respect to coordinatization of
deformation spaces. In our case, the coordinates measure how the
triangle groups are put together to form the Koebe group $\Gamma$, so
it makes sense to compute the position of $\Gamma$ in $T(\Gamma)$,
without the need of a reference fixed group. This is not difficult task,
since we have the Maskit Embedding Theorem, that reduces everything
to the one-dimensional case, and these latter groups are already
known by the previous work of this section. At this point the reader
should look at \S 2.4., where it is explained how to find subgroups
$\Gamma_j$, $1\leq j\leq 3p-3+n$ of $\Gamma$ that give the Maskit
Embedding. Given the group $\Gamma$, we decompose it in simpler
subgroups wtih one-dimensional deformation spaces, say $\Gamma_j$,
for $1\leq j\leq 3p-3+n$, and for each of these subgroups we compute
the coordinates $\alpha_j$ as explained in this section. The the
coordinate of $\Gamma$ will be $(\alpha_1,\ldots,\alpha_{3p-3+n}$.
The only point one may wonder is what happens if we start with a
different component of $\pi^{-1}(T_j)$, say $D'_j$, giving a
different decomposition of $\Gamma$. Let $\Gamma'_j$
be the corresponding stabilizer. We have that there is an element
$\gamma\in\Gamma$ such that $\gamma\Gamma_j\gamma^{-1}=\Gamma'_j$.
This means that the two stabilizers are conjugated in $PSL(2,\Bbb
C)$, and therefore the coordinate $\alpha_j$ will be the same for
both groups, since the expression of this coordinate is invariant
under conjugation by M\"{o}bius transformations (and the relation
that defines \T spaces kills conjugation by elements of $PSL(2,\Bbb C)$).

Consider now the inverse problem: let
$(\alpha_1,\ldots,\alpha_{3p-3+n})$ be a point in $T(\Gamma)$. We want
to find a Koebe group $\Gamma$ corresponding to this point. To do so,
we have to consider three different types of partition curves. Supose
we have constructed a Koebe group $\Gamma_{j-1}$, corresponding to
the first $j-1$ coordinates, $\alpha_1,\ldots,\alpha_{j-1}$.

\noindent{\bf Case $1$:} The curve $a_j$ disconnects $S$. This means
that the construction corresponding to this curve is an AFP.
Let $S_1$ and $S_2$ be the two
parts of $T_j-a_j$. We have that one of the two
parts, say $S_1$, has been already uniformized in the previous steps (i.e., in
the construction up to $a_{j-1}$).
Let $G_1=\Gamma(\mu_j,\nu(1),\nu(2);a,b,c)$ be the triangle subgroup
of $\Gamma_{j-1}$ corresponding to $S_1$.
Here $\nu(k)$ are just some ramification numbers of the
signature of $S$. Let $G_2=\Gamma(\mu_j,\nu(3),\nu(4);d,e,f)$ be a triangle
group that uniformizes orbifolds with the signature than $S_2$. All we
need to do is to find the parameters $(d,e,f)$ such that (i) $G_1$ and
$G_2$ share the element $A_j$ uniformizing $a_j$ and (ii) the
coordinate of $G_1*_{<A_j>}G_2$, as computed earlier, is precisely $\alpha_j$.
The group $\Gamma_j$ will be the group generated by $\Gamma_{j-1}$
and $G_2$.

\noindent{\bf Case $2$:} $a_j$ does not disconnect $S$, and $T_j$ is
of type $(0,4)$. The parts of $T_j$ will be uniformized by two
triangle subgroups $G_1=\Gamma(\mu_j,\nu(1),\nu(2);a,b,c)$ and
$G_2=\Gamma(\mu_j,\nu(3),\nu(4);d,e,f)$ of $\Gamma_{j-1}$. Let $A_1$
and $A_2$ be the elements of $G_1$ and $G_2$, respectively,
uniformizing $a_j$ with the correct orientations. Choose a
transformation $C$ such that $CA_2C^{-1}=A_1$, and such that the
coordinate correspoding to $G_1*_{<A_1>}CG_2C^{-1}$ is $\alpha_{j}$.
Then the group $\Gamma_j$ is the group generated by
$\Gamma_{j-1}$ and $C$.

\noindent{\bf Case $3$:} $a_j$ does not disconnects $S$ and the type of
$T_j$ is $(1,1)$. Let $G_1=\Gamma(\mu_j,\mu_j,\nu(1);a,b,c)$ be a triangle
subgroup of $\Gamma_{j-1}$ uniformizing $T_j-a_j$. Let $A$ and $B$ be
the two canonical generators of $G_1$. Find a transformation $C$ such
that $CB^{-1}C^{-1}=A$, and the coordinate of $G_1*_C$ is $\alpha_j$.
Then the group $\Gamma_j$ is generated by $\Gamma_{j-1}$ and $C$.

This algorithm completes the proof of theorems $1$ and $2$ of \S 1.
These constructions give us a way to find explicitly
generators of constructible Koebe groups. It is not hard to program a
computer to get the computations done fast and easly.

\section{Some properties of \T spaces of Koebe groups}
As we have remarked in \S 4.4, the \T space of a Koebe group is
not equivalent to the \T space of the quotient orbifold, unless all
the weights are equal to $\infty$.
So let us assume that $\Gamma$ is a constructible Koebe group
uniformizing some orbifold with a maximal partition with weights,
$(S,\cal P,\cal N)$, where at least one of the elements of $\cal N$
is finite. By results of L. Bers \cite{bers:spaces}, B. Maskit
\cite{mas:self} and I. Kra \cite{kra:spaces}, we have that the
universal covering space of $T(\Gamma)$ is the space $T(S)$.
Here we will use Maskit's
version, since his geometric description fits better in our work. See
also \cite{ear:sch} for an explanation of these results in the context
of Schottky groups (another type of Koebe groups not covered in this
paper; namely those with weights equal to $1$).

Before proceeding any further with the proof of our next result,
we need to introduce some background on surface topology.
Let $c$ be a simple closed curve on a surface $S$, and let $N_c$ be a
tubular neighborhood of $c$, homeomorphic via $h_c$, to the annulus
$A_c=\{(r,\theta);~ 1<r<3,~ 0\leq\theta\leq 2\pi\}$, with the usual
identification
modulo $2\pi$. Suppose that $c$ corresponds to the circle $r=2$. Consider
the self-homeomorphisms of $A_c$ given by
$$f_c(r,\theta)=\left\{\begin{array}{lr}(r,\theta)\, ,&1<r\leq 2\\
(r,\theta+2\pi (r-2))\, ,&2<r<3.\end{array}\right .$$
Since $f_c$ extends to the boundary of $A_c$ as the identity mapping, we can
consider the homeomorphisms of $S$ given by
$$\tau_c(r,\theta)=\left\{\begin{array}{lr}id\, ,&\mbox{on }S-A_c\\
h_c^{-1}f_ch_c\, ,&\mbox{on }N_c.\end{array}\right .$$
The Dehn twist around $c$ is just the mapping class of $\tau_c$, which we
will also denote by $\tau_c$.

Consider the minimal normal subgroup $G$ of the mapping class
group of $S$ contaning $\tau_{a_j}^{\mu_j}$ for finite $\mu_j$.
\setcounter{thm}{2}
\begin{thm}$T(\Gamma)\cong T(S)/G$\end{thm}
\setcounter{thm}{6}
\begin{pf} Let $T(S)\rightarrow T(\Gamma)$ be the universal
covering of $T(\Gamma)$.
Let $H$ denote the covering group of this mapping.
In \cite{mas:self} it is proven that $H$
consists of the elements $f$ of
the mapping class group of $S$ such that $f$ lifts to
$\tilde{f}:\Delta\rightarrow\Delta$ and
$\tilde{f}\circ\gamma=\gamma\circ\tilde{f}$, for all $\gamma\in\Gamma$.
It is well known that for each mapping class we can take $f$ to be
quasiconformal. The mapping $f$ induces a quasiconformal
self-homeomorphism on the puctured surface
$S_0=S-\{x_j;~\nu_j<\infty\}$, which we will also denote by $f$. This
mapping satisfies that for each curve $a_j$ of $\cal P$,
$f(a_j)$ is freely homotopic to some $a_{s(j)}$,
where $s$ is a permutation of $3p-3+n$ elements (this is due to the
fact that the parabolic and elliptic elements uniformizing the curves
of the partition play a special role in $\Gamma$). Then by a result
of L .Bers \cite{bers:thurs} and I. Kra (\cite{kra:horoc}
we can assume without loss of generality that $f(a_j)=a_{s(j)}$

\noindent For each $1\leq j\leq 3p-3+n$, let $T_j$ be the connected component
of $S-\{a_k;~a_k\in {\cal P},~k\neq j\}$ containing $a_j$; these are
called the modular parts of $S$. It is clear that $f$ induces a
permutation amongst the $T_j'$s.
For each $1\leq j\leq 3p-3+n$, it is not hard
to find a simple closed curve $b_j$ such that $b_j$ intersects
$a_j$ and it is disjoint from $a_k$ for any other curve of the partition
$\cal P$, $k\neq j$ (called \lq\lq dual\rq\rq~ curves in
\cite[pg. 227]{kra:variational}). If $f$ belongs to $H$, then $f$ must
commute with the elements $B_j$ uniformizing $b_j$, and therefore we have
that $f$ must preserves the sets $T_j$.

\noindent So we have reduced the problem to the one dimensional case.
Assume then that one of the modular parts, say $T_j$, is of type
$(0,4)$. The subgroup of the mapping class group of $T_j$ that
preserves the curve $a_j$ is generated by the \lq half Dehn twist\rq~
around $a_j$. Since we are assuming that the elements of
$H$ preserves all the transformations of $\Gamma$, we have that (the
restriction of $f$ to $T_j$) must be a power of $\tau_{a_j}$. The
action of this mapping on the curve $b_j$ is given by $b_j\mapsto
a_jb_ja_j^{-1}$. A similar relation holds on $\Gamma$. Therefore, if
$f\in H$, then the restriction of $f$ to $T_j$ is of the form
$\tau_{a_j}^{k\mu_j}$, with $k\in\Bbb Z$.

\noindent Similarly, for $T_j$ of type $(1,1)$, we have that $f$ reduces to
$\tau_j^k$, with $k$ integer;
but in this case the action of $f$ will be given by $b_j\mapsto a_jb_j$.
So we get $k$ is a multiple of $\mu_j$.

\noindent Since we have to conjugate in order to get the first step
of $f$ fixing
the curves $a_j$, we have the normality condition satisfied as well,
so we get $H=G$ and the theorem is proven.\end{pf}

In order to prove the last result, we need to give a
slightly different interpretation of the \T space of a Fuchsian group
$F$. See \cite{bg:isom} and \cite{ek:hol} for more details.
A {\em Beltrami coefficient} for $F$, $\mu$, is a measurable
function of norm less than $1$ such that
$(\mu\circ\gamma)\overline{\gamma'}/\gamma'=\mu$ for all
$\gamma\in F$. Given a Beltrami coefficient $\mu$, there exists a
unique quasiconformal mapping $w_\mu$ such that
$(w_\mu)_{\overline z}=\mu (w_\mu)_z$,
and $w$ fixes $\infty$, $0$ and $1$. It is not hard to see that
$w_\mu F w_\mu^{-1}$ is again a group of M\"{o}bius transformations.
Two Beltrami coefficients $\mu$ and $\nu$, are{\em equivalent} if and
only if $w_\mu=w_\nu$ on the real axis. The set of equivalence
classes of Beltrami coefficients for $F$ is the \T space of $F$,
$T(F)$. On can prove that $\mu$ and $\nu$
are equivalent if and only if there is a M\"{o}bius tranformation $A$
such that $w_\mu\circ\gamma\circ w_\nu^{-1}=
A\circ w_\nu\circ\gamma\circ w_\nu^{-1}\circ A^{-1}$, for all
$\gamma\in F$, which fits with our first definition of \T space of a
Kleinian group.

Let $S$ be an orbifold of type $(p,n)$, and let $f$ be a
quasiconformal homeomorphism of $S$ onto another orbifold $S'$, such
that $f$ respect the ramification values of the point of $S$. This
simply means that the ramification value of $x$ and $f(x)$ are equal
for all $x$ in $S$. We will call such a mapping a (quasiconformal)
{\em deformation} of $S$. Two such mappings, $f:S\rightarrow S_1$ and
$g:S\rightarrow S_2$ are equivalent if there exists a conformal
mapping $\phi:S_1\rightarrow S_2$, preserving the ramification
values, and such that $g^{-1}\circ\phi\circ f$ is homotopic to the
identity mapping on $S$. The set of equivalence classes of
quasiconformal deformations of $S$ is the \T space of $S$, $T(S)$.

Let $F$ be a Fuchsian group such that ${\Bbb H}/F\cong S$. Then
we have a natural isomorphism $\psi:T(S)\rightarrow T(F)$ given by
$f\mapsto \mu(f)$, where $f$ is a (quasiconformal) deformation of $S$
and $\mu(f)=f_{\overline z}/f_z$ is the Beltrami coefficient of $f$
computed in local coordinates.

\smallskip
\noindent{\it Proof of theorem $4$.}
Let $S$ be an orbifold of type $(p,n)$ with at least one point with
finite ramification value. Let $F$ be a Fuchsian group uniformizing $S$.
Let $S_0$ be the surface of genus $p$ with $n$
punctures obtained by removing from $S$ the points with finite
ramification value (by our assumptions we have that $S_0\neq S$). Let
$F_0$ be a Fuchsian group uniformizing $S_0$. Consider the set
${\Bbb H}_F={\Bbb H}-\{$fixed point of elliptic elemets of $F\}$.
We have that ${\Bbb H}_F/F\cong S_0$. Therefore, there
exists a universal covering map $h:{\Bbb H}\rightarrow{\Bbb H}_F$,
such that $\pi_0=\pi\circ h$, where $\pi_0:{\Bbb H}\rightarrow S_0$
and $\pi:{\Bbb H}_F\rightarrow S_0$ are the natural projection mappings.
Actually, here $\pi$ is only the
restriction of the canonical mapping from $\Bbb H$ onto $S$, but by
an abuse of notation we will denote in the same way. In \cite{bg:isom}
and \cite{ek:hol} it is proven that the mapping $h$
induces an isomorphisms between the \T spaces $T(F)$ and $T(F_0)$,
which is defined by the
formula $h^*\mu\circ h=\mu h'/\overline{h'}$, for all $\mu\in T(F_0)$.
It is not hard to understand the isomorphism $h^*$ in terms of the \T
spaces of the orbifolds $S$ and $S_0$. Let $f:S\rightarrow \tilde{S}$ be a
deformation of $S$. Restric $f$ to $S_0$, to obtain a deformation
$f|_{S_0}:=g:S_0\rightarrow\tilde{S}_0$ of $S_0$. This induces a
mapping from $r^*:T(S)\rightarrow T(S_0)$ by $f\mapsto g$. Then the mapping
$h^*$ is the inverse of $r^*$ at the group level. This means
that the following diagram is commutative:
$$\CD
T(S) \cdr{r^*}{} T(S) \\
\cdd{\psi_0}{} \cdd{}{\psi} \\
T(F_0) \cdl{}{h^*} T(F).
\endCD $$
Now the theorem follows easily from the definition of the isomorphism
$r^*$ as follows. We are given a maximal partition $\cal P$ and a set
of weights $\cal N$ on $S$. By the Maskit Existence theorem, we have
Koebe groups, $\Gamma$ and $\Gamma_0$, uniformizing
$(S,\cal P,\cal N)$ and $(S_0,\cal P,\cal N)$ in the invariant
components $\Delta$ and $\Delta_0$, respectively.
Let $\varphi:T(S)\rightarrow T(\Gamma)$ and
$\varphi_0:T(S_0)\rightarrow T(\Gamma_0)$ be the covering mapping
given by theorem $3$. The covering groups of
these mappings, $H$ and $H_0$ respectively, are the normal subgroup of
the mapping class groups of $S$ and $S_0$ generated by the Dehn twists
around the curves of $\cal P$ with finite weight. We have that $r^*$
maps $H$ onto $H_0$, by the own definition of $r^*$. So we can
project this mapping to the level of the deformation spaces of Koebe
groups, obtaining a function $\tilde{r}$, that makes the following
diagram commutative:
$$\CD
T(S) \cdr{r^*}{} T(S) \\
\cdd{\varphi_0}{} \cdd{}{\varphi} \\
T(\Gamma_0) \cdr{}{\tilde{r}} T(\Gamma).
\endCD $$
It is clear that $\tilde{r}$ is a bijection, given the desired result.
{\nolinebreak\begin{flushright}$\Box$
\end{flushright}\smallskip}

\ifx\undefined\bysame
\newcommand{\bysame}{\leavevmode\hbox to3em{\hrulefill}\,}
\fi

\sn School of Maths, Tata Institute of Fundamental Research, Bombay, India\\
pablo@motive.math.tifr.res.in

\begin{thebibliography}{10}

\bibitem{bers:spaces}
L.~Bers, {\em Spaces of {K}leinian groups}, Several {C}omplex {V}ariables
  {M}aryland 1970, Lecture {N}otes in {M}ath., vol. 155, Springer, {B}erlin,
  1970, pp.~9--34.

\bibitem{bers:thurs}
\bysame, {\em An extremal problem for quasiconformal mappings and a theorem by
  {T}hurston}, Acta Math. {\bf 141} (1978), 73--98.

\bibitem{bg:isom}
L.~Bers and L.~Greenberg, {\em Isomorphisms between {T}eichm{\"{u}}ller
  spaces}, Advances in the Theory of Riemann surfaces, Ann. of Math. Studies
  66, 1971, pp.~53--79.

\bibitem{ear:sch}
C.~Earle, {\em The group of biholomorphic self-mappings of {S}chottky space},
  Ann. Acad. Sci. Fenn. Ser. A I Math. {\bf 16} (1991), no.~2, 399--410.

\bibitem{ek:hol}
C.~Earle and I.~Kra, {\em On holomorphic mappings between {T}eichm{\"{u}}ller
  spaces}, Contributions to Analysis, Academic Press, 1974, pp.~107--124.

\bibitem{ares:thesis}
P.~Ar{\'{e}}s Gastesi, {\em Complex coordinates for the {T}eichm{\"{u}}ller
  spaces of b-groups with torsion}, Ph.D. thesis, State {U}niversity of {N}ew
  {Y}ork at {S}tony {B}rook, 1993.

\bibitem{ares:horoc1}
\bysame, {\em Coordinates for the {T}eichm{\"{u}}ller spaces of b-groups with
  torsion}, Submited to Annal Acad. Sci. Fenn. Ser. A I Math., 1993.

\bibitem{ares:bers}
\bysame, {\em The {B}ers-{G}reenberg theorem and the {M}askit embedding for
  {T}eichm{\"{u}}ller spaces}, Submited to Bull. London Math. Soc., 1994.

\bibitem{kra:spaces}
I.~Kra, {\em On spaces of {K}leinian groups}, Comment. Math. Helv. {\bf 47}
  (1972), 53--69.

\bibitem{kra:variational}
\bysame, {\em Non-variational global coordinates for {T}eichm{\"{u}}ller
  spaces}, Holomorphic functions and Moduli II, Math. Sci. Res. Inst. Publ.,
  vol.~11, Springer, 1988, pp.~221--249.

\bibitem{kra:horoc}
\bysame, {\em Horocyclic coordinates for {R}iemann surfaces and moduli spaces.
  {I}: {T}eichm{\"{u}}ller and {R}iemann spaces of {K}leinian groups}, J. Amer.
  Math. Soc. {\bf 3} (1990), 499--578.

\bibitem{km:defor}
I.~Kra and B.~Maskit, {\em The deformation space of a {K}leinian group}, Amer.
  J. Math. {\bf 103} (1981), 1065--1102.

\bibitem{mas:self}
B.~Maskit, {\em Self-maps of {K}leinian groups}, Amer. J. Math. {\bf 93}
  (1971), 840--856.

\bibitem{mas:moduli}
\bysame, {\em Moduli of marked {R}iemann surfaces}, Bull. Amer. Math. Soc. {\bf
  80} (1974), 773--777.

\bibitem{mas:class}
\bysame, {\em On the classification of {K}leinian groups: {I}-{K}oebe groups},
  Acta Math. {\bf 135} (1975), 249--270.

\bibitem{mas:kg}
\bysame, {\em Kleinian {G}roups}, Grundlehren der mathematischen
  Wissenschaften, vol. 287, Springer-Verlag, Berlin, Heidelberg, 1988.

\bibitem{mas:klein4}
\bysame, {\em On {K}lein's {C}ombination {T}heorem {IV}}, Trans. Amer. Math.
  Soc. (1993), 265--294.

\bibitem{nag:teic}
S.~Nag, {\em The {C}omplex {A}nalytic {T}heory of {T}eichm{\"{u}}ller
  {S}paces}, John Wiley \& Sons, 1988.

\bibitem{jouni:gluing}
J.~Parkkonen, {\em Elliptic {G}luing {C}onstructions}, preprint, 1994.

\end{thebibliography}
\end{document}